\documentclass[12pt,twoside,leqno]{amsart}
\usepackage{amssymb,amsbsy,amsmath,amsfonts,amssymb,amscd,epsfig,
times,graphics,color,xypic}

\definecolor{blue}{cmyk}{1.,1.,0.,0.63}
\definecolor{red}{cmyk}{0.,1.,1.,0.63}
\definecolor{green}{cmyk}{1.,0.,1.,0.63}

\usepackage[latin1]{inputenc}
\newcommand{\C}{\mathbb{C}}

\newcommand{\N}{\mathbb{N}}
\newcommand{\R}{\mathbb{R}}

\newcommand{\blue}{\textcolor{blue}}
\newcommand{\green}{\textcolor{green}}
\newcommand{\red}{\textcolor{red}}

\newtheorem{theorem}{Theorem}
\newtheorem{proposition}{Proposition}
\newtheorem{lemma}{Lemma}

\begin{document}
\parindent 0.75cm


\title[
The Hartogs extension theorem on $(n-1)$-complete complex spaces
]{
The Hartogs extension theorem 
\\
on $(n-1)$-complete complex spaces
}

\author{Jo\"el Merker and Egmont Porten}

\address{
D\'epartement de Math\'ematiques et Applications, UMR 8553
du CNRS, \'Ecole Normale
Sup\'erieure, 45 rue d'Ulm, F-75230 Paris Cedex 05, 
France. \ \
{\it Internet}:
{\tt http://www.cmi.univ-mrs.fr/$\sim$merker/index.html}}

\email{merker@dma.ens.fr}

\address{Department of Engineering, Physics and
Mathematics, Mid Sweden University, Campus Sundsvall,
S-85170 Sundsvall, Sweden
}

\email{Egmont.Porten@miun.se}

\date{\number\year-\number\month-\number\day}


\begin{abstract}
Employing Morse theory and the method of analytic discs but no
$\overline{ \partial}$ techniques, we establish a version of the
Hartogs extension theorem in a singular setting, namely: for every
domain $\Omega$ of an $(n-1)$-complete normal complex space of pure
dimension $n \geqslant 2$, and for every compact set $K \subset
\Omega$ such that $\Omega \backslash K$ is connected, holomorphic or
meromorphic functions in $\Omega \backslash K$ extend holomorphically
or meromorphically to $\Omega$.
\end{abstract}

\maketitle

\vspace{-0.5cm}

\begin{center}
\begin{minipage}[t]{11cm}
\baselineskip =0.35cm
{\scriptsize

\centerline{\bf Table of contents}

\smallskip

{\bf 1.~Introduction \dotfill 1.}

{\bf 2.~Statement of the results \dotfill 2.}

{\bf 3.~Geometrical preparations \dotfill 5.}

{\bf 4.~Holomorphic extension to $D_{\rm reg}$ \dotfill 10.}

{\bf 5.~Meromorphic extension on nonnormal complex spaces \dotfill 16.}

\smallskip

{\footnotesize\tt \hfill \blue{\bf [5 colored illustrations]}}

}\end{minipage}
\end{center}

\section*{ \S1.~Introduction}

The goal of the present article is to perform a generalization of the
classical Hartogs extension theorem in certain singular complex spaces
which enjoy appropriate convexity conditions, using the method of
analytic discs for local extensional steps and Morse-theoretical tools
for the global topological control of monodromy.

In its original form, the theorem states that in an arbitrary bounded
domain $\Omega \Subset \C^n$ ($n \geqslant 2$), every compact set $K
\subset \Omega$ with $\Omega \backslash K$ connected is an illusory
singularity for holomorphic functions, namely $\mathcal{ O}( \Omega
\backslash K) = \mathcal{ O}( \Omega ) \big\vert_{ \Omega \backslash
K}$ (for history, motivations and background, we refer {\it e.g.}  to
\cite{ hm2002, ra1986, ra2002}).  By now, the shortest proof, due to
Ehrenpreis, follows easily from the simple proposition that
$\overline{ \partial }$-cohomology with compact support vanishes in
bidegre $(0, 1)$ ({\it see} \cite{ ho1966}). Along these lines and after
results due to Kohn-Rossi, the Hartogs theorem was generalized to $(n -
1)$-complete complex manifolds by Andreotti-Hill \cite{ ah1972}, {\it
i.e.} manifolds exhausted by a $\mathcal{ C}^\infty$ function whose
Levi-form has at least $2$ positive eigenvalues at every point. We also
refer to \cite{ ks2001} for an approach via the 
holomorphic Plateau boundary problem.

To endeavor the theory in general singular complex spaces $\big( X,
\mathcal{ O}_X \big)$, it is at present advisable to look for methods
avoiding global $\overline{ \partial }$ techniques, as well as global
integral kernels, because such tools are not yet available.  The
geometric Hartogs theory was attacked long ago by Rothstein, who
introduced the notion of $q$-convexity. On the other hand, within the
modern sheaf-theoretic setting, the so-called {\sl Andreotti-Grauert
theory}\, allows to perform extension (of holomorphic functions, of
differentials forms, of coherent sheaves, {\it etc.}) from shell-like
regions of the form $\big\{ z\in X: \, a < \rho ( z) <b \big\}$ into
their inside $\big\{ z\in X : \rho ( z) < b \big\}$, where $\rho$ is a
fixed $(n-1)$-convex exhaustion function for $X$. Geometrically
speaking, one performs holomorphic extension by means of the {\sl
Grauert bump method}\, through the level sets of $\rho$ in the direction
of decreasing values, jumping finitely many times across the
critical points of $\rho$.

However, a satisfying, complete generalization of the Hartogs theorem
should apply to general excised bounded domains $\Omega \backslash K$
lying in an $(n-1)$-complete complex space $\big( X, \mathcal{ O}_X
\big)$, not only to shells $\{ a < \rho < b\}$ relative to the
$(n-1)$-convex exhaustion function. But then, after perturbing and
smoothing out $\partial \Omega$, one must unavoidably take account of
the critical points of $\rho \big\vert_{ \partial \Omega}$ and also of
the possible multi-sheetedness of the intermediate stepwise
extensions. This causes considerably more delicate topological
problems than in the well known Grauert bump method, in which
monodromy of the holomorphic (or meromorphic, or sheaf-theoretic)
extensions from $\{ a < \rho < b \}$ to $\{ a' < \rho < b \}$ with $a
' < a$ is almost freely assured\footnote{ The reader in referred to
point {\bf 2)} of the proof of Prosition~4.1 below and to Figure~3 in
Section~4 for an illustration of the concerned univalent extension
argument.}, even across critical points of $\rho$.  Considering simply
a domain $\Omega \Subset \C^n$ ($n\geqslant 2$), with obvious
exhaustion $\rho (z) := \vert \! \vert z \vert \! \vert$, the
classical Hartogs theorem based on analytic discs and on Morse theory
was worked out in \cite{ mp2006b}, where emphasis was put on rigor in
order to provide with firm grounds the subsequent works on the
subject. The essence of the present article is to transfer such an
approach to $(n - 1)$-complete general complex spaces, 
where $\overline{ \partial }$ techniques are still lacking,
with some new
difficulties due to the singularities.

\section*{\S2.~Statement of the results}

Thus, let $\big( X,\mathcal{ O}_X \big)$ be a reduced complex analytic
space of pure dimension $n \geqslant 2$, equipped with an open cover
$X = \bigcup_{ j\in J}\, U_j$ together with holomorphic isomorphisms
$\varphi_j : \ U_j \to {\sf A}_j$ onto some closed complex analytic
sets ${\sf A }_j$ contained in balls $\widetilde{ \sf B}_j \subset
\C^{ N_j}$, some $N_j \geqslant 2$. By definition (\cite{ fi1976,
gu1990}), a $\mathcal{ C}^\infty$ function $f : X \to \C$ is locally
represented as $f \vert_{ U_j} = \widetilde{ \sf f}_j \circ \varphi_j$
for some collection of $\mathcal{ C }^\infty$ ``ambient'' functions
$\widetilde{ \sf f}_j : \widetilde{ \sf B}_j \to \C$, $j\in J$. A
real-valued continuous function $\rho$ on $X$ is an {\sl exhaustion
function}\, if sublevel sets $\{ z \in X : \, \rho(z) < c\}$ are
relatively compact in $X$ for every $c\in \R$.  A $\mathcal{
C}^\infty$ function $\rho : X \to \R$ is called {\sl strongly
$q$-convex}\, if the $\mathcal{ C}^\infty$ ambient $\widetilde{
\rho}_j : \widetilde{ \sf B}_j \to \R$ can be chosen to be strongly
$q$-convex, {\it i.e.} their Levi-forms $i\, \partial \overline{
\partial} ( \widetilde{ \rho }_j)$ have at least $N_j - q + 1$
positive eigenvalues at every point, for all $j \in
J$. Finally\footnote{\, The previous definitions are known to be
independent of the choices\,\,---\,\,covering, embeddings $\varphi_j$,
dimensions $N_j$, extensions $\widetilde{ {\scriptstyle{ (\bullet )}}
}$, {\it see}\, \cite{ fi1976, gr1994, gu1990}.}, $X$ is called {\sl
$q$-complete}\, if it possesses a $\mathcal{ C }^\infty$ strongly
$q$-convex exhaustion function. Note that the $1$-complete spaces are
precisely the Stein spaces.

We will mainly work with a {\sl normal}\, $(n-1)$-complete $X$, and we
recall that a reduced complex space $\big(X, \mathcal{ O}_X \big)$ is
{\sl normal}\, if the sheaf of {\sl weakly holomorphic functions},
namely functions defined and holomorphic on 
the regular part
$X_{\rm 
reg} = X \backslash X_{\rm
sing}$ which are $L_{\rm loc}^\infty$ on $X$, coincides with the complete
sheaf $\mathcal{ O}_X$ of holomorphic functions on $X$. Then $X_{\rm
sing}$ is of codimension $\geqslant 2$ at every point of $X$ (\cite{
fi1976, gu1990}) and for every open set $U \subset X$, both
restriction maps
\def\theequation{2.1}\begin{equation}
\mathcal{O}_X(U)
\longrightarrow
\mathcal{O}_X\big(
U\backslash X_{\rm sing}
\big)
\ \ \ \ \ 
\text{\rm and}
\ \ \ \ \ 
\mathcal{M}_X(U)
\longrightarrow
\mathcal{M}_X\big(
U\backslash X_{\rm sing}
\big)
\end{equation}
are bijective\footnote{\, The first statement yields immediately that
every point $z \in X$ has a neighborhood basis $\big( \mathcal{ V}_k
\big)_{ k \in \N}$ such that $X_{\rm reg} \cap \mathcal{ V}_k$ is
connected; also, $X_{\rm reg}$ itself is connected.  The second
statement is known as Levi's extension theorem (\cite{ gere1984},
p.~185).}, where $\mathcal{ M }_X$ denotes the meromorphic sheaf. To
generalize Hartogs extension, normality of $X$ is an unavoidable
assumption, because there are examples of Stein surfaces $S$ having a
single singular point $\widehat{ p}$ which are {\it not}\, normal 
(\cite{ gu1990}, vol.~II,
p.~196),
whence $K := \{ \widehat{ p} \}$ fails to be removable for holomorphic
functions defined in a neighborhood of $K$.

We can now state our main result.

\def\thetheorem{2.2}\begin{theorem}
Let $X$ be a connected $(n-1)$-complete normal complex space of pure
dimension $n\geqslant 2$. Then for every domain $\Omega \subset X$ and
every compact set $K \subset \Omega$ with $\Omega \backslash K$
connected, holomorphic or meromorphic functions on $\Omega \backslash
K$ extend holomorphically or meromorphically and uniquely to
$\Omega${\rm :}
\[
\mathcal{O}_X(\Omega\backslash K)
=
\mathcal{O}_X(\Omega)\big\vert_{\Omega\backslash K}
\ \ \ \ \ \ \
\text{\rm or}
\ \ \ \ \ \ \
\mathcal{M}_X(\Omega\backslash K)
=
\mathcal{M}_X(\Omega)\big\vert_{\Omega\backslash K}.
\]
\end{theorem}

Some comments on the hypotheses are in order.  Firstly, connectedness
of $X$ is not a restriction, since otherwise, $\Omega$ would be
contained in a single component of $X$. Secondly, as $X$ is
$(n-1)$-complete, $i\, \partial \overline{ \partial}\, \big( \rho \big
\vert_{ X_{ \rm reg}} \big)$ has at least $2$ positive eigenvalues at
every point $z\in X_{ \rm reg}$, and
consequently, each super-level set
\[
\big\{z\in X:\, 
\rho(z)>c\big\},
\]
has a pseudoconcave boundary at every smooth point $z\in X_{\rm reg}$
with $d\rho ( z) \neq 0$ 
and in fact, 
the Levi-form of this boundary has at least one negative
eigenvalue at $z$. Thirdly, by a theorem of Ohsawa (\cite{ oh1984}), every
(connected) $n$-dimensional {\it noncompact}\, complex manifold is
$n$-complete, and in fact, easy examples show that Hartogs extension
may fail: take the product $X := R\times S$ of two Riemann surfaces,
with $R$ compact and $S$ {\it noncompact}, take a point $s \in S$ and
set $K := R \times \{ s\}$; by~\cite{ fo1981}, 
there exists a meromorphic function
function having a pole of order $1$ at $s$, whence $\mathcal{ O} ( X)$
does not extend through $K$. Consequently, in the category of strong
Levi-form assumptions, $(n-1)$-convexity is sharp.

\smallskip

For the theorem, the main strategy of proof consists of performing
holomorphic or meromorphic extension entirely within the regular part
of $X$.

\def\theproposition{2.3}\begin{proposition}
With $X$, $\Omega$ and $K$ as in Theorem~2.2,
holomorphic or meromorphic functions on $\big[ \Omega \backslash K
\big]_{ \rm reg}$ extend holomorphically or meromorphically to
$\Omega_{\rm reg}$.
\end{proposition}

Notice that both $\big[ \Omega \backslash K \big]_{\rm reg}$ and
$\Omega_{\rm reg}$ are connected (footnote~3). Then by~\thetag{ 2.1},
extension immediately holds to $\Omega$. This yields Theorem 2.2 if
one takes the proposition for granted; Sections~3 and~4 below are
devoted to prove this proposition.

\medskip

For meromorphic extension, one could in principle well avoid the
assumption of normality. In the case of meromorphic extension, we get
a general result valid for reduced spaces without further local assumptions.

\def\thetheorem{2.4}\begin{theorem}
Let $X$ be a globally irreducible $(n-1)$-complete reduced complex
space of pure dimension $n\geqslant 2$. Then for every domain $\Omega
\subset X$ and every compact set $K \subset \Omega$ with $[\Omega
\backslash K]_{\rm reg}$ connected, meromorphic functions on $\Omega
\backslash K$ extend meromorphically and uniquely to $\Omega${\rm :}
\[
\mathcal{M}_X(\Omega\backslash K)
=
\mathcal{M}_X(\Omega)\big\vert_{\Omega\backslash K}.
\]
If moreover the data lie in $\mathcal{O}_X(\Omega\backslash 
K)$, the extension is weakly holomorphic.
\end{theorem}

The proof, also relying upon an application of Proposition~2.3, is
postponed to Section~5; an example in \S5.1 shows that requiring only
that $\Omega \backslash K$ is connected does not suffices.

\smallskip

For the proposition, the main difficulty is that $X_{\rm sing}$ can in
general cross $\Omega \backslash K$. We will approach $X_{\rm sing}$
from the regular part and fill in progressively $\Omega_{\rm reg}$ by
means of the super-level sets of a suitable modification $\mu$ of the
exhaustion $\rho$, such that $\mu$ is still strongly $(n-1)$-convex
but exhausts only $X_{\rm reg}$ in a neighborhood of $\overline{
\Omega}$. To verify that the extension procedure devised in~\cite{
mp2006b} can be performed, preliminaries are
required.

\section*{\S3.~Geometrical preparations}

\subsection*{3.1.~Smoothing out the boundary}
To launch the filling procedure, we want to view the connected open
set $\Omega \backslash K$ as a neighborhood of some convenient
connected hypersurface $M$ contained in $\big( \Omega \backslash
K\big) \cap X_{\rm reg}$.

\def\thelemma{3.2}\begin{lemma}
Let $X$, $\Omega$ and $K$ be as in Theorem~2.2. Then there is a domain
$D \Subset \Omega$ containing $K$ such that $M := \partial
D \cap X_{\rm reg}$ is a $\mathcal{ C}^\infty$ {\rm connected}
hypersurface of $X_{\rm reg}$.
\end{lemma}

\proof 
Suppose first that $X = \C^n$.  Let ${\sf d}$ be a regularized
distance function (\cite{ st1970}) for $K$, {\it i.e.}  a $\mathcal{
C}^\infty$ real-valued function with $K = \{ {\sf d} = 0 \}$ and
$\frac{ 1}{ c}\, {\rm dist} \, ( x, K) \leqslant {\sf d} ( x)
\leqslant c \, {\rm dist}\, ( x, K)$ for some constant $c >1$, where
${\rm dist}$ is the Euclidean distance in $\R^{ 2n}$. By Sard's
theorem, there are arbitrarily small $\varepsilon >0$ such that
$\widehat{ M} := \{ {\sf d} = \varepsilon \}$ is a $\mathcal{
C}^\infty$ hypersurface of $\R^{ 2n}$ bounding the open set $\widehat{
\Omega} := \{ {\sf d} < \varepsilon \}$ which satisfies $K \subset
\widehat{ \Omega} \Subset \Omega$.  However, since $\widehat{ M}$ need
not be connected, we must modify it.

To this aim, we pick finitely many disjoint closed simple $\mathcal{
C}^\infty$ arcs $\gamma_1, \dots, \gamma_r$ which meet $\widehat{ M}$
transversally only at their endpoints such that $\widehat{ M} \cup
\gamma_1 \cup \cdots \cup \gamma_r$ is connected. Since $\Omega
\backslash K$ is connected, we can insure that each $\gamma_k$ is
contained in $\Omega \backslash K$.  

\begin{center}
\begin{picture}(0,0)%
\includegraphics{connecting-tubes.pstex}%
\end{picture}%
\setlength{\unitlength}{4144sp}%
\begingroup\makeatletter\ifx\SetFigFont\undefined%
\gdef\SetFigFont#1#2#3#4#5{%
  \reset@font\fontsize{#1}{#2pt}%
  \fontfamily{#3}\fontseries{#4}\fontshape{#5}%
  \selectfont}%
\fi\endgroup%
\begin{picture}(5424,1914)(439,-1537)
\put(1627,-1455){\makebox(0,0)[lb]{\smash{{\SetFigFont{9}{10.8}{\familydefault}{\mddefault}{\updefault}{\color[rgb]{0,0,0}{\bf Fig.~1: Connectifying the smoothed out boundary}}%
}}}}
\put(3767,-1222){\makebox(0,0)[lb]{\smash{{\SetFigFont{10}{12.0}{\familydefault}{\mddefault}{\updefault}{\color[rgb]{0,0,.69}\blue{$M$}}%
}}}}
\put(1202,-512){\makebox(0,0)[lb]{\smash{{\SetFigFont{8}{9.6}{\familydefault}{\mddefault}{\updefault}{\color[rgb]{.82,0,0}\red{$\gamma_1$}}%
}}}}
\put(2287,-441){\makebox(0,0)[lb]{\smash{{\SetFigFont{8}{9.6}{\familydefault}{\mddefault}{\updefault}{\color[rgb]{.82,0,0}\red{$\gamma_2$}}%
}}}}
\put(3287,-502){\makebox(0,0)[lb]{\smash{{\SetFigFont{8}{9.6}{\familydefault}{\mddefault}{\updefault}{\color[rgb]{.82,0,0}\red{$\gamma_3$}}%
}}}}
\put(4643,-619){\makebox(0,0)[lb]{\smash{{\SetFigFont{8}{9.6}{\familydefault}{\mddefault}{\updefault}{\color[rgb]{.82,0,0}\red{$\gamma_4$}}%
}}}}
\put(3931,-289){\makebox(0,0)[lb]{\smash{{\SetFigFont{9}{10.8}{\familydefault}{\mddefault}{\updefault}{\color[rgb]{0,0,0}$K$}%
}}}}
\put(1898,-863){\makebox(0,0)[lb]{\smash{{\SetFigFont{10}{12.0}{\familydefault}{\mddefault}{\updefault}{\color[rgb]{0,0,0}$K$}%
}}}}
\put(5171,-1010){\makebox(0,0)[lb]{\smash{{\SetFigFont{9}{10.8}{\familydefault}{\mddefault}{\updefault}{\color[rgb]{0,0,0}$K$}%
}}}}
\put(5108, 94){\makebox(0,0)[lb]{\smash{{\SetFigFont{10}{12.0}{\familydefault}{\mddefault}{\updefault}{\color[rgb]{0,.82,0}\green{$\widehat{M}$}}%
}}}}
\put(2976,184){\makebox(0,0)[lb]{\smash{{\SetFigFont{9}{10.8}{\familydefault}{\mddefault}{\updefault}{\color[rgb]{0,.82,0}\green{$\widehat{M}$}}%
}}}}
\end{picture}%

\end{center}

We can then modify $\widehat{ M}$
in the following way: we cut out a very small ball in $\widehat{ M}$
around each endpoint of every $\gamma_k$, and we link up the connected
components of the excised hypersurface with $r$ thin tubes $\simeq \R
\times S^{ 2n - 2}$ almost parallel to the $\gamma_k$, smoothing out
the corners appearing near the endpoints. The resulting hypersurface
$M$ is $\mathcal{ C}^\infty$ and connected. Since each $\gamma_k$ is
either contained in $\widehat{ \Omega} \cup \widehat{ M}$ or in $\R^{
2n} \big \backslash \widehat{ \Omega }$, a new open set $D$ with
$\partial D = M$ is obtained by either deleting from $\widehat{
\Omega}$ or adding to $\widehat{ \Omega}$ the thin tube around each
$\gamma_k$. All the tubes around the $\gamma_k$ which are contained in
$\R^{ 2n} \big \backslash \widehat{ \Omega}$ constitute thin open
tunnels between the components of $\widehat{ \Omega}$, whence
$D$ is connected.

\smallskip

On a general complex space $X$, the idea is to embed a
neighborhood of $\overline{ \Omega}$ smoothly into some Euclidean space
$\R^N$ and then to proceed similarly.

We can assume that the holomorphic isomorphisms $\phi_j : U_j \to {\sf
A}_j \subset \widetilde{ \sf B}_j \subset \C^{ N_j}$ are defined in
slightly larger open sets $U_j' \Supset U_j$, for all $j \in J$.  Pick
$\mathcal{ C}^\infty$ functions $\lambda_j$ having compact support in
$U_j'$ and satisfying $\lambda_j = 1$ on $\overline{
U}_j$; prolong them to be $0$
on $X$ outside $U_j'$. By compactness, there is a finite open cover:
\[
\overline{\Omega}
\subset
U_{j_1}\cup\cdots\cup 
U_{j_m}.
\]
Consider the $\mathcal{ C}^\infty$ map, valued in $\R^N$ with $N:= 2(
N_{j_1} + \cdots + N_{ j_m}) + m$, which is defined by:
\[
\Psi
:=
\big(
\lambda_{j_1}\cdot\phi_{j_1},\dots,
\lambda_{j_m}\cdot\phi_{j_m},\ \
\lambda_{j_1},\dots,\lambda_{j_m}
\big).
\]
It is an immersion at every point $x$ of $U_{ j_1} \cup \cdots \cup
U_{ j_m}$, because $x$ belongs to some $U_{ j_k}$, whence the $j_k$-th
component $\lambda_{ j_k} \cdot \phi_{ j_k} \equiv \phi_{ j_k}$ of
$\Psi$ is even an embedding of $U_k \ni x$. Furthermore, we claim that
$\Psi$ separates points. Indeed, if we set:
\[
W_{j_k}
:=
\big\{
z\in X:\ 
\lambda_{j_k}(z)=1
\big\},
\]
then clearly $U_{ j_k} \subset W_{ j_k} \subset U_{ j_k}'$.  Pick two
distinct points $x, y \in U_{ j_1} \cup \cdots \cup U_{ j_m}$. Then
$x$ belongs to some $U_{ j_k}$, so $\lambda_{ j_k} ( x) = 1$. If
$\lambda_{ j_k} ( y) \neq 1$, then $\Psi( y) \neq \Psi ( x)$ and we
are done. If $\lambda_{ j_k} (y) = 1$, {\it i.e.}  if $y \in W_{
j_k}$, then the $j_k$-th component of $\Psi$ distinguishes $x$ from
$y$, since $\lambda_{j_k} \cdot \phi_{j_k} (y) = \phi_{ j_k} ( y)$
differs from $\phi_{ j_k} (x)$ because $\phi_{ j_k}$ embeds $U_{ j_k}
'$ into $\R^{ 2 N_{ j_k}}$. So $\Psi$ embeds into $\R^N$ the
neighborhood $U_{j_1}
\cup \cdots \cup U_{ j_m}$ of $\overline{ \Omega}$.

\smallskip

We choose a regularized distance function ${\sf d}_{\Psi (K)}$ for
$\Psi ( K)$ in $\R^N$. We stratify $X$ so that $X_{\rm reg}$ is the
single largest stratum (remind it is connected) and then stratify
$X_{\rm sing}$ by listing all connected components of $\big[ X_{\rm
sing} \big]_{\rm reg}$, then continuing with $\big[ X_{\rm sing}
\big]_{\rm sing}$, and so on inductively.  By Sard's theorem and the
stratified transversality theorem (\cite{ hi1976}), for almost every
$\varepsilon >0$, the intersection
\[
\big\{ x \in \R^N : \ {\sf
d}_{\Psi(K)} (x) = \varepsilon \big\} \cap 
\Psi\big(\Omega_{\rm reg}\big)
\]
is a $\mathcal{ C }^\infty$ real hypersurface of $\Psi ( \Omega_{\rm
reg})$ having finitely many connected components which are contained
in $\Psi \big( [ \Omega \backslash K ]_{\rm reg} \big)$. Importantly,
we can construct the thin connecting tubes so that they {\it lie all
entirely inside}\, $\Psi \big( \big[ \Omega \backslash K \big]_{\rm
reg} \big)$, thanks to the fact that $\Psi \big( \Omega_{\rm reg}
\big)$ is locally (arcwise) connected, also near points of $\Psi \big(
\Omega_{\rm sing} \big)$. Then the remaining arguments are the same
and we put everything back to $X$ {\it via}\, $\Psi^{ - 1}$, getting a
connected $\mathcal{ C}^\infty$ hypersurface $M \subset [\Omega
\backslash K]_{\rm reg}$ and a domain $D$ with $K \subset D \Subset
\Omega$.  (We remark that normality of $X$ was crucially used.)
\endproof

As we said, we will perform the filling procedure entirely inside
$X_{\rm reg}$. This is possible thanks to an idea of Demailly which
consists of modifying the initial exhaustion $\rho$ so that $X_{\rm
sing}$ is put at $- \infty$. A recent application 
of this idea also appears in~\cite{ df2007}.

\subsection*{3.3.~Putting $X_{\rm sing}$ into a well}
By Lemma~5 in~\cite{ de1990}, there exists an {\sl almost
plurisubharmonic function}\footnote{\, {\it i.e.} by definition, a
function which is locally the sum of a psh function and of a
$\mathcal{ C}^\infty$ function, or equivalently, a function $v$ whose
complex Hessian $i \, \partial \overline{ \partial}\, v$ has bounded
negative part.} $v$ on $X$ which is $\mathcal{ C}^\infty$ on $X_{\rm
reg}$ and has poles along $X_{\rm sing}$:
\[
X_{\rm sing}
=
\big\{
v=-\infty
\big\}.
\]
As in Section~2, if ${\sf A}_j = \varphi_j (U_j)$ is represented in a
local ball $\widetilde{ \sf B}_j \subset \C^{ N_j}$ of radius $r_j>0$
centered at $z_j \in \C^{ N_j}$ as the zero-set $\{ g_{j,\nu} = 0\}$
of finitely many functions $g_{ j,\nu}$ holomorphic in a neighborhood
of the closure of $\widetilde{ \sf B}_j$, the local ambient
$\widetilde{ v}_j : \widetilde{ \sf B}_j \to \{ - \infty \} \cup \R$
is essentially of the form\footnote{\, In addition, a regularized
maximum function (\cite{ de1990}) is used to smoothly glue these
different definitions on all finite intersections $A_{ j_1} \cap
\cdots \cap A_{ j_m}$ and the formula given here is exact on a
sub-ball $\widetilde{ \sf C}_j \subset \widetilde{\sf B}_j$.}:
\[
\widetilde{ v}_j
=
{\rm log}\, 
\Big(
\sum_\nu\,\vert 
g_{j,\nu} 
\vert^2 
\Big)
-
\frac{1}{r_j^2-\vert z-z_j\vert^2}.
\]
Thus, locally
on each $\widetilde{ \sf B}_j$, the function $v$ we pick
from~\cite{ de1990} is of the form:
\[
\widetilde{v}_j
=
\widetilde{u}_j
+
\widetilde{\sf r}_j,
\]
with $\widetilde{ u}_j$ strictly psh, $\mathcal{ C}^\infty$ on
$\widetilde{ \sf B}_j \big \backslash \big[ {\sf A}_j \big]_{\rm
sing}$, equal to $\{ - \infty \}$ on $\big[ {\sf A}_j \big]_{\rm
sing}$ and with a remainder $\widetilde {\sf r}_j$ which is $\mathcal{
C}^\infty$ on the whole of $\widetilde{ \sf B}_j$.  Notice that each
$\widetilde{ v}_j$ is $L_{\rm loc}^\infty$.

\subsection*{3.4.~Modified strongly 
$(n-1)$-convex exhaustion function $\mu$} Pick a constant $C > 0$ such
that $\max_{ \overline{ D}}\, (\rho) < C$.

\def\thelemma{3.5}\begin{lemma}
There exists $\varepsilon_0 >0$ such that for all $\varepsilon$ with
$0 < \varepsilon \leqslant \varepsilon_0$, the function{\rm :}
\[
\mu
:=
\rho+\varepsilon\,v
\]
is $\mathcal{ C}^\infty$ on $X_{\rm reg}$ and satisfies{\rm :}

\begin{itemize}

\smallskip\item[{\bf (a)}]
$\max_{ \overline{ D} } \, (\mu) < C${\rm ;}

\smallskip\item[{\bf (b)}]
$X_{\rm sing} = \{ \mu = - \infty \}${\rm ;}

\smallskip\item[{\bf (c)}]
$\mu$ is strongly $(n-1)$-convex in a neighborhood of $\{ \rho
\leqslant C\}$.

\end{itemize}\smallskip
\end{lemma}

\proof
Property {\bf (b)} holds provided only that $\varepsilon < \frac{ C -
\max_{ \overline{ D}} \, (\rho) }{\max_{ \overline{ D}}\,
( v) }$. Furthermore, {\bf (a)} is clear since $\rho$ is $\mathcal{
C}^\infty$ and since $X_{\rm sing} = \{ v = - \infty \}$.  To check
{\bf (c)}, we compute Levi-forms as $(1, 1)$-forms:
\def\theequation{3.6}\begin{equation}
\aligned
i\,\partial\overline{\partial}\,
\widetilde{\mu}_j
&
=
i\,\partial\overline{\partial}\,
\widetilde{\rho}_j
+
\varepsilon\,
i\,\partial\overline{\partial}\,
\widetilde{v}_j
\\
&
=
i\,\partial\overline{\partial}\,
\widetilde{\rho}_j
+
\varepsilon\,
i\,\partial\overline{\partial}\,
\widetilde{u}_j
+
\varepsilon\,
i\,\partial\overline{\partial}\,
\widetilde{\sf r}_j.
\endaligned
\end{equation}
Here, $\varepsilon\, i\, \partial \overline{ \partial }\, \widetilde{
u}_j$ adds positivity to $i\, \partial \overline{ \partial }\,
\widetilde{ \rho}_j$ (since $\widetilde{ u}_j$ is psh), whereas the
negative contribution due to $i\, \partial \overline{ \partial}\,
\widetilde {\sf r}_j$ is bounded from below on $\{ \rho \leqslant 2C
\}$, and consequently, $\varepsilon >0$ can be chosen small enough so
that $i\, \partial \overline{ \partial }\, \widetilde{ \mu}_j$ still
has $2$ eigenvalues $>0$ at every point.
\endproof

In the next section, while applying the holomorphic extension
procedure of~\cite{ mp2006b}, we shall have to insure that the
extensional domains attached to $M$ from either the outside or the
inside {\it cannot go beyond $\{ \rho \leqslant C\}$}.  So we have to
prepare in advance the curvature of the limit hypersurface $\{ \rho =
C\} \cap X_{\rm reg}$.

Enlarging $C$ of an arbitrarily small increment if necessary, we can
assume (thanks to Sard's theorem) that $C$ is a regular value of $\rho
\big \vert_{ X_{\rm reg}}$, so that
\[
\Lambda
:=
\{\rho=C\}
\cap 
X_{\rm reg}
\] 
is a $\mathcal{ C }^\infty$ real hypersurface of $X_{\rm reg}$.

\def\thelemma{3.7}\begin{lemma}
Lowering again $\varepsilon >0$ if necessary, the
following holds{\rm :}

\begin{itemize}

\smallskip\item[{\bf (d)}]
{\it At every point $q$ of the $\mathcal{ C}^\infty$ real hypersurface
$\Lambda = \{ \rho = C \} \cap X_{\rm reg}$, one can find a complex
line $E_q \subset T_q^c \Lambda$ on which the Levi-forms of both
$\rho$ and $\mu$ are positive}.

\end{itemize}\smallskip

\end{lemma}

Here, $q \mapsto E_q$ might well be discontinuous, but this shall not
cause any trouble in the sequel.

\proof
Each $p \in \{ \rho = C\}$ is contained in some $U_{ j ( p)}$, whence
$\rho$ is represented by an ambient function $\widetilde{ \rho}_{ j (
p)} : \widetilde{ \sf B}_{ j (p)} \to \R$ whose Levi-form has at least
$N_{ j ( p)} - n + 2$ eigenvalues $>0$. By diagonalizing the Levi
matrix $i\, \partial \overline{ \partial} \widetilde{ \rho}_{ j(p)}$
at the central point of $\widetilde{ \sf B}_{ j (p)}$, we may easily
define, in some small open sub-ball $\widetilde{ \sf C}_{ j (p)}
\subset \widetilde{ \sf B}_{ j (p)}$ having the same center, a
$\mathcal{ C}^\infty$ family $\widetilde{ q} \mapsto \widetilde{ F}_{
\widetilde{ q}}$ of complex $(N_{j(p)}-n+2)$-dimensional
affine subspaces such that the Levi-form
of $\widetilde{ \rho}_{ j(p)}$ is positive definite on every
$\widetilde{ F}_{ \widetilde{ q}}$, for every $\widetilde{ q} \in
\widetilde{ \sf C}_{ j (p)}$. 

Next, if we set $V_{ j(p)} := \varphi_{ j(p)}^{ - 1} \big(
\widetilde{\sf C}_{ j(p)} \big)$, which is an open subset of $U_{
j(p)}$, we can cover the compact set $\{ \rho = C\}$ by finitely many
$V_{ j(p)}$, hence there is a finite number of points $p_a$, $a = 1,
\dots, A$, such that
\[
\{\rho=C\}
\subset
V_{j(p_1)}\cup\cdots\cup
V_{j(p_A)}.
\]
According to~\thetag{ 3.6}, on each 
$\widetilde{ \sf C}_{ j( p_a)}$, $a = 1, \dots A$, 
we have:
\[
i\,\partial\overline{\partial}\,
\widetilde{\mu}_{j(p_a)}
=
i\,\partial\overline{\partial}\,
\widetilde{\rho}_{j(p_a)}
+
\varepsilon\,i\,\partial\overline{\partial}\,
\widetilde{u}_{j(p_a)}
+
\varepsilon\,i\,\partial\overline{\partial}\,
\widetilde{\sf r}_{j(p_a)}.
\]
We choose $\varepsilon >0$ so small that the remainder $\varepsilon
\,i \, \partial \overline{ \partial}\, \widetilde{\sf r}_{ j(p_a)}$
does not perturb positivity on 
$\widetilde{ \sf C}_{ j ( p_a)}$ for every $a = 1, \dots A$, and we get
that $i\, \partial \overline{
\partial}\, \widetilde{\mu }_{ j(p_a)}$ is still positive on
$\widetilde{ F}_{ \widetilde{ q}}$ for every $\widetilde{ q} \in
\widetilde{ \sf C}_{ j ( p_a )}$, and every $a = 1, \dots A$.

Let $q \in \{ \rho = C\} \cap 
X_{\rm reg}$. Then $q \in V_{ j ( p_a)}$ for some $a$. 
We set $\widetilde{ q} := \varphi_{ j( p_a)} ( q)\in \widetilde{ \sf
C}_{ j ( p_a)}$ and we define:
\[
F_q
:=
\big(
d\,\varphi_{j(p_a)}
\big)^{-1}
\Big(
\widetilde{F}_{\widetilde{q}}\cap
T_{\widetilde{q}}\,{\sf A}_{j(p_a)}
\Big).
\] 
Then the complex linear spaces $\widetilde{
F}_{ \widetilde{ q}}$ and $F_q$ are at least of dimension $2$ and
the Levi-form of $\mu$ is positive on any $1$-dimensional subspace
$E_q \subset F_q \cap T_q^c \, \Lambda$.

\endproof

Next, applying Morse transversality theory, we may perturb $\mu$ in
$X_{\rm reg} \cap \{ \rho < 2\,C\}$ in an arbitrarily small way, so
that\footnote{\, The previous four properties being preserved,
especially {\bf (d)} on $\{ \rho = C\}$.}:

\begin{itemize}

\smallskip\item[{\bf (e)}]
{\it $\mu$ is a Morse function on $X_{\rm reg} \cap \{ \rho < 2\,C\}$
having finitely many or at most countably many critical points{\rm ;}
moreover, different critical points of $\mu$ are located in
different level sets $\{ \mu = c \}$.}

\end{itemize}\smallskip

Of course, if they are infinite in number, critical values can only
accumulate at $- \infty$. Similarly, we may perturb $\rho$ very
slightly near $\{ \rho = C\}$ so that:

\begin{itemize} 

\smallskip\item[{\bf (f)}]
the $\mathcal{ C}^\infty$ hypersurface $\{ \rho = C \} \cap X_{\rm
reg}$ does not contain any critical point of $\mu$.

\end{itemize}\smallskip

Finally, again thanks to Morse transversality theory, we may perturb
the connected $\mathcal{ C}^\infty$ hypersurface $M \subset \partial
D$ of Lemma~3.2 
in an arbitrarily small way so that\footnote{\, The perturbed $M$
being still contained in $\{ \rho < C \}$ and in the original
connected corona $\Omega \backslash K$.}:

\begin{itemize}

\smallskip\item[{\bf (g)}]
{\it $M$ does not contain critical points of $\mu$, and $\mu \big
\vert_M$ is a Morse function on $M$ having finitely many or at most
countably many critical points{\rm ;} moreover, any two different
critical points of $\mu$ or of $\mu \big \vert_M$ have different
critical values.}

\end{itemize}\smallskip

We draw a diagram, where $X_{\rm sing}$ is symbolically represented as
a continuous broken line having spikes, with a level-set $\{ \mu =
\widehat{ c}\}$ which is critical for $\mu \big \vert_M$ and a single
critical point $\widehat{ p} \in M \cap \{ \mu = \widehat{ c}\}$.

\begin{center}
\begin{picture}(0,0)%
\includegraphics{mu-M-X-sing.pstex}%
\end{picture}%
\setlength{\unitlength}{4144sp}%
\begingroup\makeatletter\ifx\SetFigFont\undefined%
\gdef\SetFigFont#1#2#3#4#5{%
  \reset@font\fontsize{#1}{#2pt}%
  \fontfamily{#3}\fontseries{#4}\fontshape{#5}%
  \selectfont}%
\fi\endgroup%
\begin{picture}(5424,2049)(439,-1648)
\put(5311,-1128){\makebox(0,0)[lb]{\smash{{\SetFigFont{10}{12.0}{\familydefault}{\mddefault}{\updefault}{\color[rgb]{.82,0,0}\red{$X_{\rm sing}$}}%
}}}}
\put(5306, 59){\makebox(0,0)[lb]{\smash{{\SetFigFont{10}{12.0}{\familydefault}{\mddefault}{\updefault}{\color[rgb]{.82,0,0}\red{$X_{\rm sing}$}}%
}}}}
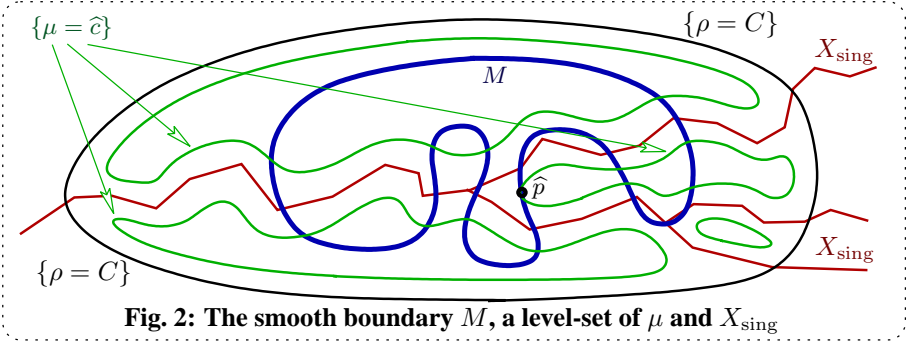
\put(1178,-1554){\makebox(0,0)[lb]{\smash{{\SetFigFont{10}{12.0}{\familydefault}{\mddefault}{\updefault}{\color[rgb]{0,0,0}{\bf Fig.~2: The smooth boundary $M$, a level-set of $\mu$ and $X_{\rm sing}$}}%
}}}}
\put(3318,-98){\makebox(0,0)[lb]{\smash{{\SetFigFont{9}{10.8}{\familydefault}{\mddefault}{\updefault}{\color[rgb]{0,0,.69}\blue{$M$}}%
}}}}
\put(3621,-764){\makebox(0,0)[lb]{\smash{{\SetFigFont{9}{10.8}{\familydefault}{\mddefault}{\updefault}{\color[rgb]{0,0,0}$\widehat{p}$}%
}}}}
\put(617,196){\makebox(0,0)[lb]{\smash{{\SetFigFont{9}{10.8}{\familydefault}{\mddefault}{\updefault}{\color[rgb]{0,.82,0}\green{$\{\mu=\widehat{c}\}$}}%
}}}}
\put(4517,214){\makebox(0,0)[lb]{\smash{{\SetFigFont{10}{12.0}{\familydefault}{\mddefault}{\updefault}{\color[rgb]{0,0,0}$\{\rho = C\}$}%
}}}}
\put(649,-1278){\makebox(0,0)[lb]{\smash{{\SetFigFont{10}{12.0}{\familydefault}{\mddefault}{\updefault}{\color[rgb]{0,0,0}$\{\rho = C\}$}%
}}}}
\end{picture}%

\end{center}

\section*{\S4.~Holomorphic extension to $D_{\rm reg}$}

For $c\in \R$, we introduce
\[
X_{\mu>c}:= 
\{z\in X:\,
\mu(z)
>
c\}.
\] 
This open set is contained in $X_{\rm reg}$, since $X_{\rm sing} = \{
\mu = - \infty\}$. For every connected component $M_{\mu > c}'$ of
\[M_{\mu>c}
:=
M\cap
X_{\mu>c}= 
M\cap
\{\mu>c\}, 
\]
we want to fill in (by means of a finite number of families of
analytic discs) a certain domain $Q_{ \mu > c}'$ which is enclosed by
$M_{ \mu > c}'$ inside $\{ \mu > c\}$.  Similarly as in
Proposition~5.3 of~\cite{ mp2006b}, we must consider {\it all}\, the
connected components $M_{\mu > c}'$ and analyze the combinatorics of
how they merge or disappear.

Let $\mathcal{ V} ( M)$ be a thin tubular neighborhood of $M$, whose
thinness shrinks to zero while approaching $X_{\rm sing}$.  For every
connected component $M_{ \mu > c}'$ of $M_{ \mu > c}$, we denote by
$\mathcal{ V} \big( M_{\mu > c}' \big)_{ \mu > c}$ the part of
$\mathcal{ V} ( M)$ around $M_{ \mu > c}'$ again intersected with $\{
\mu > c\}$. It is a connected tubular neighborhood of $M_{ \mu > c}'$
inside $\{ \mu > c\}$.

\def\theproposition{4.1}\begin{proposition}
Let $c \in \R$ with $c < \max_M\, (\mu) < C$ be any regular value of
$\mu$ and of $\mu \big \vert_M$. Let $M_{ \mu > c}'$ be any nonempty
connected component of $M \cap X_{\mu > c}$. Then there is a unique
connected component $Q_{ \mu > c}'$ of $X_{\mu > c} \big \backslash
M_{ \mu > c}'$ which is relatively compact in $X_{\rm reg}$ and
contained in $\{ \rho < C\}$ with the property that two different
domains $Q_{\mu > c}'$ and $Q_{\mu > c}''$ are either disjoint or one
is contained in the other. Furthermore, for every holomorphic or
meromorphic function $f$ defined in the thin tubular neighborhood
$\mathcal{ V} ( M)$ of $M$, there exists a unique holomorphic or
meromorphic extension $F$, constructed by means of a finite number of
$(n-1)$-concave Levi-Hartogs figures and defined in
\[
Q_{\mu>c}'\bigcup
\mathcal{V}\big(M_{\mu>c}'\big)_{\mu>c},
\]
such that $F = f$ when both functions
are restricted to $\mathcal{V}\big(M_{\mu>c}'\big)_{\mu>c}$.
\end{proposition}

\proof
We only describe the modifications one must bring to the arguments
of~\cite{ mp2006b}. 

\medskip

{\bf 1)}
The Levi-form of the compact $\mathcal{ C}^\infty$ boundary $\{ \mu =
c \}$ of the super-level set $\{ \mu > c\}$ (contained in $X_{\rm
reg}$) has $1$ negative eigenvalue, so that the Levi extension theorem
with analytic discs ({\it cf.}  the survey \cite{ mp2006a}) applies
at each point of $\{ \mu = c\}$. In Section~3 of~\cite{ mp2006b}, we
defined $(n-a)$-concave Hartogs figures for $1 \leqslant a \leqslant
n-1$, but we used only $1$-concave ones, because the Levi-form of
exterior of spheres $\{ \vert \!  \vert z \vert \! \vert < r\}$ in
$\C^n$ had $(n-1)$ negative eigenvalues.  Here, we start from
$(n-1)$-concave Hartogs figures, we modify them similarly as in
Section~3 of~\cite{ mp2006b} (details are skipped) and we call them
{\sl $(n-1)$-concave Levi-Hartogs figures}.

Next, we use a finite number of these figures, {\it via}\, some local
charts of $X_{\rm reg}$, to cover $\{\mu = c\}$ and to show that
holomorphic\footnote{\, Since the configuration is always local and
biholomorphic to $\C^n$ ($n = \dim X_{\rm reg}$) and since holomorphic
envelopes coincide with meromorphic envelopes in $\C^n$, meromorphic
functions enjoy exactly the same extension properties.  Thus, in
\cite{ mp2006b}, results stated for holomorphic functions are
immediately true for meromorphic functions too.}  (or meromorphic)
functions in $\{ \mu > c\}$ extend to a slightly deeper super-level
set $\{ \mu > c - \eta\}$ (provided no critical point of $\mu$ or of
$\mu \big \vert_M$ is encountered in the shell $\{ c 
\geqslant \mu > c - \eta \}$),
for some $\eta > 0$ which depends on $X$, on $n$, on $\mu$, but not on
$c$.

\medskip

{\bf 2)} Contrary to the $\C^n$ case treated in~\cite{ mp2006b}, $\mu$
may have critical points on $X_{\rm reg}$. Grauert's theory shows how
to jump across them with $\overline{ \partial}$ techniques, and we
summarize how we can proceed here\footnote{\, We emphasize
that, from the point of view of holomorphic extension, jumping across
critical points of $\mu$ on $X_{\rm reg}$ is much simpler than jumping
across critical points of $\mu \big\vert_M$, {\it cf.} the $\C^n$
case \cite{ mp2006b}.}, using only analytic discs in Levi-Hartogs figures.

Consider a point $\widehat{ p} \in X_{\rm reg}$ which is critical:
$d\mu ( \widehat{ p}) = 0$, and set $\widehat{c} := \mu ( \widehat{
p})$.  The Morse lemma provides local real coordinates centered at
$\widehat{ p}$ in which $\mu = x_1^2 + \cdots + x_k^2 - y_1^2 - \cdots
- y_{ 2n - k}^2$, for some $k$. Since $i\, \partial \overline{
\partial}\, \mu$ has at least $2$ positive eigenvalues everywhere, $k$
is $\geqslant 2$. This is a crucial fact, because this implies that
super-level sets $\{ \mu > \widehat{ c} + \delta \}$ are all
connected\footnote{\, In $\R^3$ already, this is true for the
``exterior'' $x^2 + y^2 - z^2 > \delta$ of the standard cone.} in a
neighborhood of $\widehat{ p}$, for every $\delta \in \R$ close to $0$,
and moreover, that these domains grow regularly and continuously as
$\delta$ decreases from positive values to negative values.

\begin{center}
\begin{picture}(0,0)%
\includegraphics{double-excised.pstex}%
\end{picture}%
\setlength{\unitlength}{4144sp}%
\begingroup\makeatletter\ifx\SetFigFont\undefined%
\gdef\SetFigFont#1#2#3#4#5{%
  \reset@font\fontsize{#1}{#2pt}%
  \fontfamily{#3}\fontseries{#4}\fontshape{#5}%
  \selectfont}%
\fi\endgroup%
\begin{picture}(5424,2184)(439,-1783)
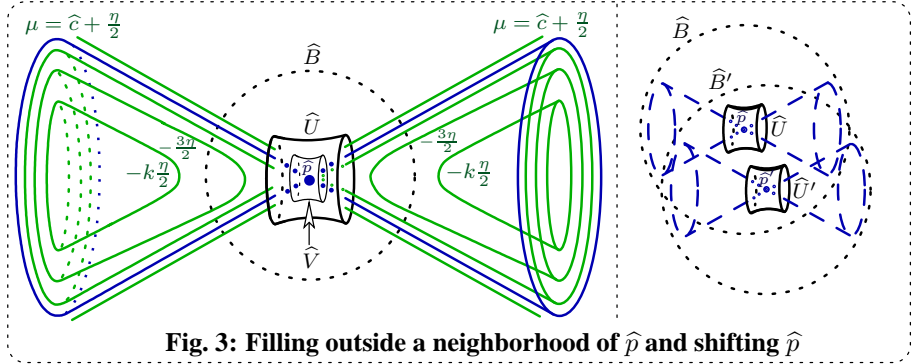
\put(1396,-1700){\makebox(0,0)[lb]{\smash{{\SetFigFont{10}{12.0}{\familydefault}{\mddefault}{\updefault}{\color[rgb]{0,0,0}{\bf Fig.~3: Filling outside a neighborhood of $\widehat{ p}$ and shifting $\widehat{p}$}}%
}}}}
\put(2215,-1203){\makebox(0,0)[lb]{\smash{{\SetFigFont{8}{9.6}{\familydefault}{\mddefault}{\updefault}{\color[rgb]{0,0,0}$\widehat{V}$}%
}}}}
\put(2218, 16){\makebox(0,0)[lb]{\smash{{\SetFigFont{8}{9.6}{\familydefault}{\mddefault}{\updefault}{\color[rgb]{0,0,0}$\widehat{B}$}%
}}}}
\put(2212,-631){\makebox(0,0)[lb]{\smash{{\SetFigFont{6}{7.2}{\familydefault}{\mddefault}{\updefault}{\color[rgb]{0,0,.69}\blue{$\widehat{p}$}}%
}}}}
\put(2220,-400){\makebox(0,0)[lb]{\smash{{\SetFigFont{8}{9.6}{\familydefault}{\mddefault}{\updefault}{\color[rgb]{0,0,0}$\widehat{U}$}%
}}}}
\put(3356,227){\makebox(0,0)[lb]{\smash{{\SetFigFont{8}{9.6}{\familydefault}{\mddefault}{\updefault}{\color[rgb]{0,.82,0}\green{$\mu=\widehat{c}+\frac{\eta}{2}$}}%
}}}}
\put(2913,-473){\makebox(0,0)[lb]{\smash{{\SetFigFont{6}{7.2}{\familydefault}{\mddefault}{\updefault}{\color[rgb]{0,.82,0}\green{$-\frac{3\eta}{2}$}}%
}}}}
\put(1348,-513){\makebox(0,0)[lb]{\smash{{\SetFigFont{6}{7.2}{\familydefault}{\mddefault}{\updefault}{\color[rgb]{0,.82,0}\green{$-\frac{3\eta}{2}$}}%
}}}}
\put(3072,-665){\makebox(0,0)[lb]{\smash{{\SetFigFont{8}{9.6}{\familydefault}{\mddefault}{\updefault}{\color[rgb]{0,.82,0}\green{$-k\frac{\eta}{2}$}}%
}}}}
\put(4810,-333){\makebox(0,0)[lb]{\smash{{\SetFigFont{6}{7.2}{\familydefault}{\mddefault}{\updefault}{\color[rgb]{0,0,.69}\blue{$\widehat{p}$}}%
}}}}
\put(5013,-416){\makebox(0,0)[lb]{\smash{{\SetFigFont{8}{9.6}{\familydefault}{\mddefault}{\updefault}{\color[rgb]{0,0,0}$\widehat{U}$}%
}}}}
\put(4424,173){\makebox(0,0)[lb]{\smash{{\SetFigFont{8}{9.6}{\familydefault}{\mddefault}{\updefault}{\color[rgb]{0,0,0}$\widehat{B}$}%
}}}}
\put(4642,-135){\makebox(0,0)[lb]{\smash{{\SetFigFont{8}{9.6}{\familydefault}{\mddefault}{\updefault}{\color[rgb]{0,0,0}$\widehat{B}'$}%
}}}}
\put(4943,-700){\makebox(0,0)[lb]{\smash{{\SetFigFont{6}{7.2}{\familydefault}{\mddefault}{\updefault}{\color[rgb]{0,0,.69}\blue{$\widehat{p}'$}}%
}}}}
\put(5147,-787){\makebox(0,0)[lb]{\smash{{\SetFigFont{8}{9.6}{\familydefault}{\mddefault}{\updefault}{\color[rgb]{0,0,0}$\widehat{U}'$}%
}}}}
\put(556,219){\makebox(0,0)[lb]{\smash{{\SetFigFont{8}{9.6}{\familydefault}{\mddefault}{\updefault}{\color[rgb]{0,.82,0}\green{$\mu=\widehat{c}+\frac{\eta}{2}$}}%
}}}}
\put(1154,-671){\makebox(0,0)[lb]{\smash{{\SetFigFont{8}{9.6}{\familydefault}{\mddefault}{\updefault}{\color[rgb]{0,.82,0}\green{$-k\frac{\eta}{2}$}}%
}}}}
\end{picture}%

\end{center}

Next, we fix a ball $\widehat{ B}$ centered at $\widehat{ p}$ and
we cut out a small neighborhood $\widehat{ U}
\subset \widehat{ B}$ of $\widehat{ p}$. 
If $\widehat{ V} \subset \widehat{ U}$ is
a small neighborhood, 
we consider the $\mathcal{ C}^\infty$ hypersurface:
\[
\big\{
\mu 
>
\widehat{c}
+ 
{\textstyle 
\frac{\eta}{2}}\}\,
\big\backslash\,
\widehat{V}.
\]
Placing finitely many $(n-1)$-concave Levi-Hartogs figures at points
of this hypersurface, we get holomorphic or meromorphic extension to
$\big\{ \mu > \widehat{ c} - \frac{ \eta}{ 2} \big\}\, \big
\backslash\, \widehat{ V}_1$, where $\widehat{ V}_1 \subset \widehat{
V}$ is slightly bigger than $\widehat{ V}$. Repeating the filling
process finitely many times until $\big\{ \mu = \widehat{ c} - \frac{
k \eta}{ 2} \big\}$ does not intersect $\widehat{ B}$, where $k$ is an
odd integer, we fill in $\widehat{ B} \big \backslash \widehat{ U}$.
At each step, monodromy of the extension is assured thanks to
connectedness of $\big\{ \mu > \widehat{ c} + \delta \big\}\, \big
\backslash\, \widehat{U}$, for every small $\delta \in \R$.  However,
we cannot fill in $\widehat{ U}$ directly this way.

The trick is to shift $\widehat{ p}$. One introduces a $\mathcal{
C}^\infty$ perturbation $\mu'$ of $\mu$ localized near $\widehat{ p}$
(namely $\mu ' = \mu$ elsewhere) such that $\mu'$ has another critical
point $\widehat{ p}'$ (having the same Morse index of course), with
corresponding neighborhoods disjoint: $\widehat{ U} \cap \widehat{ U}'
= \emptyset$ and both contained in $\widehat{ B} \cap \widehat{ B}'$.
We repeat the Levi-Hartogs filling with $\mu'$, getting holomorphic or
meromorphic extension $\big\{ \mu' > \widehat{ c} - k' \frac{ \eta}{
2} \big\}\, \big \backslash\, \widehat{ U}'$, a domain which contains
$\widehat{ B}' \big \backslash \widehat{ U}'$, hence contains
$\widehat{ U}$. Monodromy is again well controlled, just because
topologically, $\widehat{ B} \big \backslash \widehat{ U}$ and
$\widehat{ B}' \big \backslash \widehat{ U}'$ are complete shells.

\medskip

{\bf 3)} We prove the proposition by decreasing $c$.  Provided $c$
does not cross critical values of $\mu \big \vert_M$, the domains $Q_{
\mu > c}'$ do grow regularly and continuously, even when $c$ crosses
critical values of $\mu$, according to what has been said just
above. At a critical value $\widehat{ c}$ of $\mu \big \vert_M$, for a
domain $Q_{\mu > \widehat{ c}}$ whose closure contains the
corresponding unique critical point $\widehat{ p} \in M$, similarly as
in~\cite{ mp2006b}, three cases may occur:

\begin{itemize}

\smallskip\item[{\bf (i)}]
the domain $Q_{\mu > \widehat{ c} + \delta}'$ grows regularly and
continuously as $\delta$ decreases in a neighborhood of $0${\rm ;}

\smallskip\item[{\bf (ii)}]
precisely when $\delta$ becomes negative, the domain $Q_{\mu >
\widehat{ c} + \delta}'$ is merged with a second domain $Q_{\mu >
\widehat{ c} + \delta}''$ whose closure also contains $\widehat{ p}$
for $\delta = 0$ (the case of three domains or more never occurs){\rm ;}

\smallskip\item[{\bf (iii)}]
the domain $Q_{\mu > \widehat{ c} + \delta}'$ is contained in a bigger
domain $Q_{\mu > \widehat{ c} + \delta}''$ for all small $\delta > 0$,
and exactly at $\delta = 0$, the closure of the domain $Q_{\mu >
\widehat{ c}}'$ is subtracted from $Q_{\mu > \widehat{ c}}''$,
yielding a new domain $Q_{\mu > \widehat{ c}}'''$ which starts to grow
regularly and continuously as $Q_{\mu > \widehat{ c} +
\delta}'''$ for small $\delta < 0$.

\end{itemize}\smallskip

We then check by decreasing induction on $c$ that such domains are
relatively compact and are either disjoint or one is contained in the
other, and we achieve extension by means of $(n-1)$-concave
Levi-Hartogs figures similarly as in~\cite{ mp2006b}.  But here, a
single fact remains to be established, namely that the domains $Q_{\mu
> c}'$ remain all contained inside the relatively compact region $\{
\rho < C\}$.

This is true at the beginning of the filling process, namely for $c$
slightly smaller than $\max_M \, (\mu)$, because $M_{ \mu >c}$ is then
diffeomorphic to a small spherical cap (hence connected) and the
relatively compact domain enclosed by $M_{ \mu >c}$ in $X_{\mu > c}
\big \backslash M_{ \mu > c}$ is just the piece $D_{\mu > c}$ of $D$,
which is diffeomorphic to a thin cut out piece of ball close to $M$
and clearly contained in $\{ \rho < C\}$, since $ D\cup M \subset \{
\rho < C\}$ by {\bf (a)}.

To prove that all $Q_{ \mu > c}'$ are contained in $\{ \rho < C\}$, we
proceed by contradiction. Let $c^*$ be first $c$ (as $c$ decreases)
for which some $Q_{\mu > c}'$ is not contained in $\{ \rho < C\}$.  In
the process described above of constructing the domains $Q_{ \mu >
c}'$, the only discontinuity occurs in {\bf (iii)} and it consists of
a suppression.  Consequently, the domains $Q_{\mu > c}'$ cannot jump
discontinuously across $\{ \rho = C\}$, hence at $c = c^*$ (which
might be either critical or noncritical), all $Q_{\mu > c^*}'$ are
still contained in $\{ \rho \leqslant C\}$ and the boundary of at
least one domain, say $Q_{ \mu > c^*}^*$, touches the $\mathcal{
C}^\infty$ border hypersurface $\{ \rho = C\} \cap X_{\rm reg}$.

\begin{center}
\begin{picture}(0,0)%
\includegraphics{touch-mu-rho.pstex}%
\end{picture}%
\setlength{\unitlength}{4144sp}%
\begingroup\makeatletter\ifx\SetFigFont\undefined%
\gdef\SetFigFont#1#2#3#4#5{%
  \reset@font\fontsize{#1}{#2pt}%
  \fontfamily{#3}\fontseries{#4}\fontshape{#5}%
  \selectfont}%
\fi\endgroup%
\begin{picture}(5608,1824)(439,-1423)
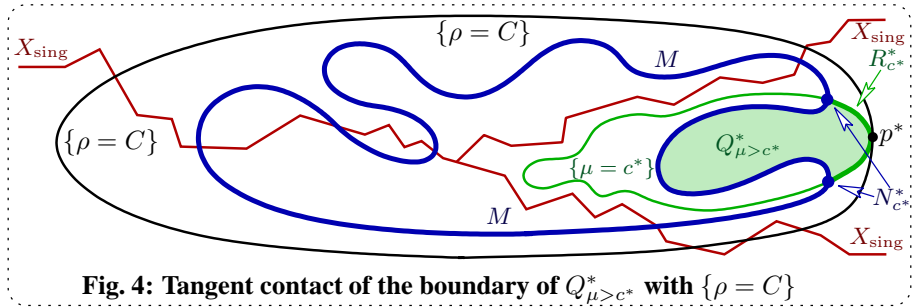
\put(893,-1339){\makebox(0,0)[lb]{\smash{{\SetFigFont{10}{12.0}{\familydefault}{\mddefault}{\updefault}{\color[rgb]{0,0,0}{\bf Fig.~4: Tangent contact of the boundary of $Q_{\mu>c^*}^*$ with   $\{\rho=C\}$}}%
}}}}
\put(3005,177){\makebox(0,0)[lb]{\smash{{\SetFigFont{10}{12.0}{\familydefault}{\mddefault}{\updefault}{\color[rgb]{0,0,0}$\{\rho=C\}$}%
}}}}
\put(777,-475){\makebox(0,0)[lb]{\smash{{\SetFigFont{10}{12.0}{\familydefault}{\mddefault}{\updefault}{\color[rgb]{0,0,0}$\{\rho=C\}$}%
}}}}
\put(484, 89){\makebox(0,0)[lb]{\smash{{\SetFigFont{9}{10.8}{\familydefault}{\mddefault}{\updefault}{\color[rgb]{.82,0,0}\red{$X_{\rm sing}$}}%
}}}}
\put(3308,-924){\makebox(0,0)[lb]{\smash{{\SetFigFont{9}{10.8}{\familydefault}{\mddefault}{\updefault}{\color[rgb]{0,0,.69}\blue{$M$}}%
}}}}
\put(4313, 17){\makebox(0,0)[lb]{\smash{{\SetFigFont{9}{10.8}{\familydefault}{\mddefault}{\updefault}{\color[rgb]{0,0,.69}\blue{$M$}}%
}}}}
\put(5663,-429){\makebox(0,0)[lb]{\smash{{\SetFigFont{10}{12.0}{\familydefault}{\mddefault}{\updefault}{\color[rgb]{0,0,0}$p^*$}%
}}}}
\put(5482,-1045){\makebox(0,0)[lb]{\smash{{\SetFigFont{9}{10.8}{\familydefault}{\mddefault}{\updefault}{\color[rgb]{.82,0,0}\red{$X_{\rm sing}$}}%
}}}}
\put(3803,-621){\makebox(0,0)[lb]{\smash{{\SetFigFont{8}{9.6}{\familydefault}{\mddefault}{\updefault}{\color[rgb]{0,.82,0}\green{$\{\mu=c^*\}$}}%
}}}}
\put(5621,-813){\makebox(0,0)[lb]{\smash{{\SetFigFont{9}{10.8}{\familydefault}{\mddefault}{\updefault}{\color[rgb]{0,0,.69}\blue{$N_{c^*}^*$}}%
}}}}
\put(4680,-486){\makebox(0,0)[lb]{\smash{{\SetFigFont{9}{10.8}{\familydefault}{\mddefault}{\updefault}{\color[rgb]{0,.82,0}\green{$Q_{\mu>c^*}^*$}}%
}}}}
\put(5594, 25){\makebox(0,0)[lb]{\smash{{\SetFigFont{9}{10.8}{\familydefault}{\mddefault}{\updefault}{\color[rgb]{0,.82,0}\green{$R_{c^*}^*$}}%
}}}}
\put(5491,182){\makebox(0,0)[lb]{\smash{{\SetFigFont{9}{10.8}{\familydefault}{\mddefault}{\updefault}{\color[rgb]{.82,0,0}\red{$X_{\rm sing}$}}%
}}}}
\end{picture}%

\end{center}

On the other hand, by definition and by construction, for each $c$,
the boundary of each $Q_{ \mu> c}'$ consists of two parts: $M_{ \mu >
c}'$, which is contained in $M$, hence remains always in $\{ \rho <
C\}$, together with a certain closed region $R_{ \mu = c}' \cup N_{\mu
= c}'$ contained in $\{ \mu = c\}$, with $R_{\mu = c}'$ open and
$N_{\mu = c}'$ being the boundary in $\{ \mu = c\}$ of $R_{\mu =
c}'$. In fact, similarly as in Section~5 of~\cite{ mp2006b}, $R_{\mu =
c}'$ is always contained in $\{ \mu = c \} \big \backslash M$ and
$N_{\mu = c}'$, always contained in $M \cap \{ \mu = c \}$ is a
$\mathcal{ C}^\infty$ real submanifold of $X_{\rm reg}$ of codimension
$2$ provided $c$ is noncritical for $\mu \big\vert_M$, while $N_{\mu =
c}'$ may have as a single singular (corner) point $\widehat{ p}$ for $c =
\widehat{ p}$ critical. But since $N_{\mu = c}'$ is a subset of $M \cap
\{ \mu = c\}$, it is always contained in $\{ \rho < C\}$.

Consequently, the boundary of $Q_{\mu > c^*}^*$ can touch $\{ \rho =
C\}$ only at some point $p^* \in R_{ \mu = c^*}^*$. So we have $\mu (
p^*) = c^*$ and $\rho (p^*) = C$, namely $p^*$ lies in $\{ \mu =
c^*\}$ and on the $\mathcal{ C}^\infty$ hypersurface $\{ \rho = C\}$.

By {\bf (f)} above, $p^* \in \{ \rho = C\}$ cannot be a critical point
of $\mu$, whence $\{ \mu = c^* \}$ and $\{ \rho = C\}$ are both
$\mathcal{ C}^\infty$ real hypersurfaces passing through $p^*$.
Furthermore, $\{ \mu \geqslant c^*\}$ is still contained in $\{ \rho
\leqslant C\}$, by definition of $c^*$, whence $T_{p^*} \{ \rho = C\}
= T_{ p^*} \{ \mu = c^*\}$.

Thanks to {\bf (d)}, there is a complex line
\[
E_{p^*} 
\subset 
T_{p^*}^c\{\rho=C\}
=
T_{p^*}^c\{\mu=c^*\}
\]
on which the Levi-forms of both $\rho$ and $\mu$ are positive
definite. On the other hand, since $\{- \mu < - c^* \}$ is contained
in $\{ \rho < C \}$, the Levi-form of $- \mu$ in the direction of $E_{
p^*}$ should then be $\geqslant$ the Levi-form of $\rho$ in the same
direction.  This is a contradiction, and the proof 
that all $Q_{\mu > c}'$ remain in $\{ \rho < C\}$
is completed. This finishes our proof of
Proposition~4.1.
\endproof

\subsection*{4.2.~End of proof of Proposition~2.3} As in Section~2
of~\cite{ mp2006b}, one checks that extension holds from $\big[ \Omega
\backslash K \big]_{\rm reg}$ to $\Omega_{\rm reg}$ provided
holomorphic or meromorphic functions defined in the thin tubular
neighborhood $\mathcal{ V} ( M)$ of $M \subset X_{\rm reg}$ do extend
uniquely to $D_{\rm reg} \bigcup \mathcal{ V} ( M)$.  So we work with
$M$, $\mathcal{ V} ( M)$ and $D_{\rm reg}$, and since everything is
exhausted as $c \to -\infty$, the conclusion of the proof of
Proposition~2.3 is an immediate consequence of the following.

\def\theproposition{4.3}\begin{proposition}
For every regular value $c > - \infty$ of $\mu \big\vert_M$,
holomorphic or meromorphic functions defined in $\mathcal{ V} ( M)$ do
extend holomorphically or meromorphically and uniquely to
\[
D_{\mu>c}
\bigcup
\mathcal{V}\big(M_{\mu>c}\big)_{\mu>c}.
\]
\end{proposition}

\proof
We set $c_1 := \max_M ( \mu) = \max_{ \overline{ D}} (\mu) < C$. There
is a unique ``$\mu$-farthest point'' $p_1 \in M$ with $\mu (p_1) =
c_1$ and this point is obviously a critical point of Morse index equal
to $- (2n-1)$ for $\mu \big\vert_M$, by virtue of {\bf (g)}.
Consequently, for all $c < c_1$ close to $c_1$, there is a single
connected component in $M_{ \mu > c}$, namely $M_{ \mu > c}$ itself,
which is diffeomorphic to a small spherical cap and encloses the
domain $D_{ \mu > c}$, diffeomorphic to a thin cut out piece of ball.
For such $c < c_1$ close to $c_1$, the proposition is thus a direct
consequence of the previous Proposition~4.1.

For arbitrary noncritical $c$, there is a well defined connected
component $M_{\mu > c}^1$ of $M_{\mu > c}$ with $p_1 \in M_{\mu
>c}^1$, and we denote by $M_{\mu > c}^2, \dots, M_{\mu > c}^k$ the
other connected components of $M_{\mu > c}$.  Also, each connected
component $D_{\mu > c}^\sim$ of $D_{\mu > c}$ is bounded by some of
the $M_{\mu > c}^j$, inside $\{ \mu > c\}$.
The problem is that the various extensions provided
by Proposition~4.1 need not stick together, but fortunately, 
we can go to deeper super-level sets $\{ \mu > c'\}$. 

\def\thelemma{4.4}\begin{lemma}
For every $c'$ with $- \infty < c' \leqslant c$ which is noncritical
for $\mu \big\vert_M$, the $\mu$-farthest point $p_1$ belongs to a
unique connected component $M_{\mu > c'}'$ of $M \cap \{ \mu > c' \}$
and the enclosed domain $Q_{\mu > c'}'$ constructed by Proposition~4.1
contains $D$ in a neighborhood of $p_1$.
\end{lemma}

\proof
Indeed, if this were not true, there would exist the first $c' = c^*$
(as $c' \leqslant c$ decreases) for which $Q_{\mu > c'}'$ switches to
the other side of $M$ near $p_1$.  According to the topological
combinatorial processus {\bf (i)}, {\bf (ii)}, {\bf (iii)} above, this
could only occur in case {\bf (iii)} with $c^*$ critical, where a
component is suppressed from a bigger one $Q_{\mu > c^*}''$ bounded by
some $M_{\mu > c^*}''$, the suppressed component necessarily being
$Q_{\mu > c^*}'$ itself. Then the bigger component $Q_{\mu > c^*}''$
would contain the side of $M$ which is exterior to $D$ near $\widehat{
p}_1$, whence
\[
c_1''
:= 
\max
\big\{
\mu(q):\
q\in 
M_{\mu>c^*}''
\big\}
\]
would necessarily be $> c_1$, which contradicts $c_1 = \max_M \,
(\mu)$.
\endproof

\begin{center}
\begin{picture}(0,0)%
\includegraphics{deepening.pstex}%
\end{picture}%
\setlength{\unitlength}{4144sp}%
\begingroup\makeatletter\ifx\SetFigFont\undefined%
\gdef\SetFigFont#1#2#3#4#5{%
  \reset@font\fontsize{#1}{#2pt}%
  \fontfamily{#3}\fontseries{#4}\fontshape{#5}%
  \selectfont}%
\fi\endgroup%
\begin{picture}(5424,2274)(439,-1873)
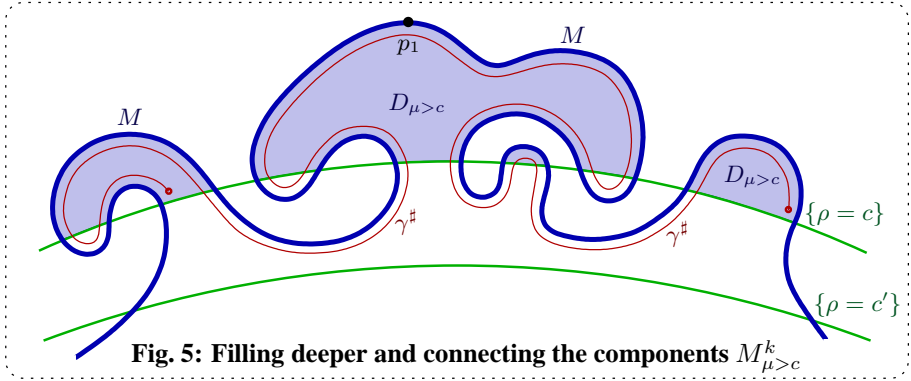
\put(1203,-1767){\makebox(0,0)[lb]{\smash{{\SetFigFont{10}{12.0}{\familydefault}{\mddefault}{\updefault}{\color[rgb]{0,0,0}{\bf Fig.~5: Filling deeper and connecting the components $M_{\mu >c}^k$}}%
}}}}
\put(2740,-255){\makebox(0,0)[lb]{\smash{{\SetFigFont{9}{10.8}{\familydefault}{\mddefault}{\updefault}{\color[rgb]{0,0,.69}\blue{$D_{\mu>c}$}}%
}}}}
\put(2794,105){\makebox(0,0)[lb]{\smash{{\SetFigFont{9}{10.8}{\familydefault}{\mddefault}{\updefault}{\color[rgb]{0,0,0}$p_1$}%
}}}}
\put(4746,-687){\makebox(0,0)[lb]{\smash{{\SetFigFont{9}{10.8}{\familydefault}{\mddefault}{\updefault}{\color[rgb]{0,0,.69}\blue{$D_{\mu>c}$}}%
}}}}
\put(4401,-1047){\makebox(0,0)[lb]{\smash{{\SetFigFont{9}{10.8}{\familydefault}{\mddefault}{\updefault}{\color[rgb]{.82,0,0}\red{$\gamma^\sharp$}}%
}}}}
\put(2778,-972){\makebox(0,0)[lb]{\smash{{\SetFigFont{9}{10.8}{\familydefault}{\mddefault}{\updefault}{\color[rgb]{.82,0,0}\red{$\gamma^\sharp$}}%
}}}}
\put(1116,-343){\makebox(0,0)[lb]{\smash{{\SetFigFont{9}{10.8}{\familydefault}{\mddefault}{\updefault}{\color[rgb]{0,0,.69}\blue{$M$}}%
}}}}
\put(3766,147){\makebox(0,0)[lb]{\smash{{\SetFigFont{9}{10.8}{\familydefault}{\mddefault}{\updefault}{\color[rgb]{0,0,.69}\blue{$M$}}%
}}}}
\put(5231,-916){\makebox(0,0)[lb]{\smash{{\SetFigFont{9}{10.8}{\familydefault}{\mddefault}{\updefault}{\color[rgb]{0,.82,0}\green{$\{\rho=c\}$}}%
}}}}
\put(5286,-1461){\makebox(0,0)[lb]{\smash{{\SetFigFont{9}{10.8}{\familydefault}{\mddefault}{\updefault}{\color[rgb]{0,.82,0}\green{$\{\rho=c'\}$}}%
}}}}
\end{picture}%

\end{center}

Next, since $M$ is connected (according to 
Lemma~3.2), we can pick a $\mathcal{
C}^\infty$ Jordan arc $\gamma$ running in $M$ which starts at $p_1$
and visits every other connected component $M_{\mu > c}^2, \dots,
M_{\mu > c}^k$ of $M_{\mu > c}$.  Since $\gamma$ is compact, there is
a noncritical $c' > - \infty$ such that $\gamma \subset \{ \mu >
c'\}$.  Fix such a $c'$ and denote by $M_{\mu > c'}'$
the connected component of $M \cap \{ \mu > c'\}$ to 
which $p_1$ belongs. Then let 
$Q_{\mu > c'}'$ be as
in Lemma~4.4.

\def\thelemma{4.5}\begin{lemma}
The domain $Q_{\mu > c'}'$ contains $D_{\mu > c}$.
\end{lemma}

\proof
Near $p_1$, this domain already contains a piece of $D$ thanks to
Lemma~4.4. From the beginning, $M$ is oriented, since it bounds the
domain $D$.  Thus, we can push $\gamma$ slightly inside $D$, getting a
curve $\gamma^\sharp$ almost parallel to $\gamma$ which is entirely
contained in $D$, and also contained in $\{ \mu > c'\}$ if the push is
sufficiently small.  Furthermore, $\gamma^\sharp$ is also entirely
contained in $Q_{\mu>c'}'$, because the extensional domain $Q_{\mu >
c'}'$ is, at least near $p_1$, located on the same side (with respect
to $M$) as $D$.

Let $D_{\mu > c}^\sim$ be any connected component of $D_{\mu > c}$. By
construction, $\gamma^\sharp$ visits $D_{\mu > c}^\sim$. Thus, every
point of $D_{\mu > c}^\sim$ may be joined to some point of
$\gamma^\sharp$ by means of some auxiliary $\mathcal{ C}^\infty$ curve
running in $D_{\mu > c}^\sim$. All such auxiliary curves do not meet
$M$, hence they do not meet $M_{\mu > c'}'$, whence they all run in
$Q_{\mu > c'}'$.  Consequently, by means of $\gamma^\sharp$ and of the
auxiliary curves in each $D_{\mu > c}^\sim$, we may connect, without
crossing $M$ even once, every point of $D_{\mu > c}$ with the starting
point of $\gamma^\sharp$, contained in $Q_{\mu > c'}'$ near
$p_1$. Thus $D_{\mu > c}$ is effectively contained in $Q_{\mu > c'}'$.
\endproof

To conclude, an application of Proposition~4.1 yields unique extension
to $Q_{\mu > c'}' \bigcup \mathcal{ V} \big( M_{\mu > c'} ' \big)_{
\mu > c'}$, and by plain restriction, we get unique extension to
$D_{\mu > c} \bigcup \mathcal{ V} \big( M_{\mu > c} \big)_{\mu > c}$.

\smallskip

This completes the proofs of Propositions~4.3 and~2.3.
\endproof

\section*{\S5.~Meromorphic extension on nonnormal complex spaces}

\subsection*{5.1.~An example}
To see that the weaker assumption that $\Omega \backslash K$ is
connected does not suffice, we consider $X = \C^2/ \big( (-1,0) \sim
(+1,0) \big)$, the euclidean $\C^2$ with two points identified. If we
define the structure sheaf by $\mathcal{ O }_{ \C^2, z}$ at all single
points and by $\mathcal{ O }_{\C^2, \pm } = \big\{ (f, g) \in
\mathcal{ O}_{ \C^2, -1} \times \mathcal{ O }_{ \C^2,1 } : f(-1,0) =
g( +1, 0) \big \}$ at the double point $(\pm 1,0)$, the space $\big(
X, \mathcal{ O }_X \big)$ is reduced and modelled near $(\pm 1,0)$ on
$\{( z,w) \in \C^2 \times \C^2 :z =w\}$.  This makes it easy to check
that the function $|z_1+1|^2 + |z_1-1|^2 + |z_2 |^2$ descends to a
1-convex exhaustion of $X$ via the quotient projection
$\pi:\C^2\rightarrow X$. Letting $\Omega := X$ and $K := \pi \big(\{
|z_1+1|^2 + |z_2|^2 = 1\} \big)$, we see that $\Omega \backslash K$ is
connected. Furthermore, $\mathcal{ O }( \Omega \backslash K)$ consists
of all functions holomorphic in $\C^2 \big \backslash \big\{ |z_1 +
1|^2 + |z_2|^2 = 1 \big\}$ which satisfy $f (-1,0) = f(+1,0)$. Then
obviously, the conclusion of Theorem 2.4 does not hold.

\subsection*{5.2.~Proof of Theorem~2.4}
To begin with, we observe that Proposition 2.3 carries over without
change to the more general setting of Theorem 2.4: indeed, thanks to
the connectedness of $[ \Omega \backslash K ]_{\rm reg}$, we may
construct $M$ and $D$ as in Lemma~3.2; the construction of an
almost psh function $v$ with $X_{\rm sing} = \{ v = - \infty \}$ holds
without assumption of normality (\cite{ de1990}), and then
Propositions~4.1 and~4.3 do go through (notice that both $\Omega
\backslash K$ and $\Omega_{\rm reg}$ are connected). Thus $\mathcal{
M}_X( \Omega \backslash K)$ extends uniquely as $\mathcal{ M }_X \big(
\Omega_{ \rm reg} \cup \left[ \Omega \backslash K \right] \big)$,
holomorphicity being preserved.

Extension across $\Omega_{ \rm sing} \cap K$ is slightly more
complicated than in the normal case due to the fact that $\Omega_{ \rm
sing }$ may have components of codimension one.  Let $\pi: \widehat{
X} \rightarrow X$ be the normalization of $X$.  Let $X_{ \rm norm }$
be the set of the normal points of $X$. Recall that $\pi$ restricts to
a biholomorphism on $\pi^{-1} (X_{ \rm norm})$. Topologically, $\pi$
is a local homeomorphism over irreducible points of $X$ and separates
the irreducible local components at reducible points.  For every open
$U\subset X$, setting $\widehat{ U } = \pi^{ -1} (U)$, we have a
canonical isomorphism $\pi^*: \mathcal{M}_X (U) \rightarrow \mathcal{
M }_{ \widehat{ X } }( \widehat{ U})$ (\cite{ gere1984},
p.~155). Hence it is enough to extend from $\mathcal{ M }_{ \widehat{
X }} \big( \widehat{ \Omega}\backslash L \big)$ to $\mathcal{ M }_{
\widehat{ X }}( \widehat{ \Omega})$, where $\widehat{ \Omega } :=
\pi^{-1}( \Omega )$ and $L := \pi^{ -1} \big( \Omega_{ \rm sing} \cap
K \big)$.

By the Levi extension theorem, we can extend through all points of
$z\in L$ with $\dim_z \pi^{ -1}( \Omega_{ \rm sing })\leqslant n-2$.
Let $H$ be an irreducible component of $\Omega_{ \rm sing }$ of
codimension one. Since $\dim \widehat{ \Omega }_{ \rm sing}\leqslant
n-2$, it follows that $\widehat{ H }' := \pi^{ -1}(H) \cap
\widehat{\Omega }_{ \rm reg }$ is dense, open and connected in
$\widehat{ H} = \pi^{-1} (H)$. Because $X$ is $(n-1)$-convex, it cannot
contain any compact analytic hypersurface according to Lemma~5.3 just
below, and $H$ has to intersect $\Omega\backslash K$. For dimensional
reasons, $\widehat{ H }'$ intersects $\big[ \pi^{-1}( \Omega \backslash
K) \big]_{ \rm reg }$, and we can apply the following version of the
Levi extension theorem for complex manifolds (\cite{ gh1978}): {\it
Let Y be an analytic subset of a complex manifold of $M$ of
codimension at least one. If $U\subset M$ is a domain containing $M
\backslash Y$ and intersecting each irreducible one-codimensional
component of $Y$, then holo-\big/meromorphic functions on $U$ extend
holo-\big/meromorphically to $M$.}

The remaining part of the singularity lies in $\widehat{ \Omega}_{ \rm
sing }$ and can be removed by the Levi extension theorem. If the
original function on $\Omega \backslash K$ is holomorphic, the
extension on $\widehat{ \Omega}$ is so too, and its push-forward to
$\Omega$ is weakly holomorphic. The proof of Theorem 2.4 is complete.
\qed

\def\thelemma{5.3}\begin{lemma}
An $(n-1)$-convex complex space $X$ of pure dimension $n$
cannot contain any analytic hypersurface $Y$ which is compact.
\end{lemma}

\proof
Let $\rho$ be an $(n-1)$-convex exhaustion function. Let $\big( U_j
\big)_{ j\in J}$ be a locally finite covering of $X$ by open subsets
which can be embedded onto analytic subsets ${\sf A}_j$ of euclidean
domains $\widetilde{ \sf B}_j \subset \C^{ N_j}$ such that the
push-forward of $\rho$ extends as an $(n-1)$-convex function
$\widetilde{ \rho}_j \in \mathcal{ C}^\infty \big ( \widetilde{ \sf
B}_j \big)$.  By an inductive deformation of $\rho$, we may arrange
that all $\widetilde{ \rho}_j$ can be chosen to be Morse functions
without critical points on ${\sf A}_j$.

If there is a compact analytic hypersurface $Y \subset X$, then
$\rho|_Y$ attains a global maximum at some point $z_0 \in Y$.  We can
assume that $z_0$ lies in some ball $\widetilde{\sf B}_j$, we denote
by ${\sf E}_j \subset {\sf A}_j \subset \widetilde{\sf B}_j \subset
\C^{ N_j}$ the local representative of $Y$ and we drop the index $j$,
because the rest of the argument is local. By construction $\big\{ z:
\, \widetilde{ \rho} (z) = \widetilde{ \rho} (z_0) \big\}$ is a smooth
$(n-1)$-convex real hypersurface such that ${\sf E } \subset\{
\widetilde{ \rho} \leqslant \widetilde{ \rho} (z_0) \}$.  Bending this
hypersurface a little, we can arrange that ${\sf E}$ is in fact
contained in $\{ \widetilde{ \rho} < \widetilde{ \rho} (z_0) \} \cup
\{ z_0 \}$ near $z_0$.  By $(n-1)$-convexity of $\widetilde{ \rho}$,
there is a piece $\Lambda$ of a small $(N - n + 1 )$-dimensional
complex plane passing through $z_0$ and contained in the complex
tangent plane $T_{ z_0 }^c \{ \widetilde{ \rho} = \widetilde{ \rho} (
z_0) \}$ on which the Levi-form $i \, \partial \overline{ \partial}
\widetilde{ \rho}$ is positive. Thus $\Lambda$ is contained in $\{
\widetilde{ \rho} > \widetilde{ \rho} ( z_0)\} \cup \{ z_0 \}$ and has
a contact of order exactly two with $\{ \widetilde{ \rho} =
\widetilde{ \rho} ( z_0) \}$ at $z_0$.  Furthermore, if we pick a
nonzero vector $v\in T_{z_0}\C^N$ which points into $\{ \rho >
\rho(z_0) \}$, the translates $\Lambda_\epsilon := \Lambda +
\varepsilon \, v$ do all lie in $\{\rho > \rho (z_0)\}$ for every
small $\varepsilon >0$, whence $\Lambda_\varepsilon \cap {\sf E}$ is
empty. But given that $\Lambda_0 \cap Y = \{z_0\} \neq\emptyset$, this
contradicts the persistence, under perturbation, of the intersection
of two complex analytic sets of complementary dimensions in
$\C^N$. The lemma is proved.
\endproof

\vfill\end{document}